\input amstex
\documentstyle{amsppt}
\language0
\NoBlackBoxes
\def\Tr{\operatorname{Tr}}
\let\ii\i
\let\ccdot\cdot
\def\cdot{\hbox to 2.5pt{\hss$\ccdot$\hss}}
\define\tsum{\tsize\sum}
\define\newquad{\hbox{\hphantom{${}={}$}}}
\define\>{\rightarrow}
\define\<{\leftarrow}
\define\[{\lbrack}
\define\]{\rbrack}
\redefine\o{\circ}

\define\al{\alpha}
\define\be{\beta}
\define\ga{\gamma}
\define\de{\delta}

\define\ze{\zeta}
\define\et{\eta}
\define\th{\theta}

\define\ka{\kappa}

\define\si{\sigma}
\define\ta{\tau}
\define\ph{\varphi}

\define\ps{\psi}
\define\om{\omega}
\define\Ga{\Gamma}

\define\La{\Lambda}

\define\Ph{\Phi}

\define\Om{\Omega}

\define\ad{{\operatorname{ad}}}
\define\Ad{{\operatorname{Ad}}}
\let\ii\i
\redefine\i{^{-1}}
\define\row#1#2#3{#1_{#2},\ldots,#1_{#3}}
\define\rowup#1#2#3{#1^{#2},\ldots,#1^{#3}}
\define\irow#1#2#3{{#1_{#2}\ldots#1_{#3}}}
\define\irowup#1#2#3{{#1^{#2}\ldots#1^{#3}}}
\define\x{\times}
\define\nmb#1#2{#2}      
\def\idx{}               


\redefine\L{{\Cal L}}
\predefine\SS\S
\redefine\S{{\Cal S}}

\define\dd#1{\tfrac \partial{\partial #1}}
\define\field#1#2#3{{#1}^{#2}\frac\partial{\partial{#1}_{#3}}}

\define\pol#1 #2#3{L^{#1}_{\text{sym}}(\Bbb R^#2,\Bbb R^#3)}

\def\today{\ifcase\month\or
 January\or February\or March\or April\or May\or June\or
 July\or August\or September\or October\or November\or December\fi
 \space\number\day, \number\year}
\topmatter
\title Invariant operators on manifolds with almost Hermitian
symmetric structures, \\II. Normal Cartan connections\endtitle
\author Andreas \v Cap, Jan Slov\'ak, Vladim\'\ii r Sou\v cek\endauthor
\abstract
In the first part of this series of papers we developed the invariant
differentiation with respect to a Cartan connection, we described
this procedure in the terms of the underlying principal
connections, and we discussed applications of
this theory to the construction of natural operators. In this part we will
extend the results of \cite{Ochiai, 70} on the existence and the
uniqueness of the so called normal Cartan connections on
manifolds with almost Hermitian symmetric structures to
first order structures which do not admit a torsion free linear
connection. Moreover, for each of these structures we obtain explicit
(universal) formulae for these
canonical connections in the terms of the curvatures of the
underlying linear connections.
\endabstract
\address
Institut f\"ur Mathematik, Universit\"at Wien, Strudlhofgasse 4,
1090 Wien, Austria\newline\indent
Department of Algebra and Geometry, Masaryk University in Brno,
Jan\'a\v ckovo n\' am\. 2a, 662~95 Brno, Czech
Republic\newline\indent
Mathematical Institute, Charles University, Sokolovsk\'a 83,
Praha, Czech Republic
\endaddress
\thanks This work was mostly done during the stay of the authors at the Erwin
Schr\"odinger Institute in Vienna.
The second author is also supported by the GA\v CR, grant Nr\.
201/93/2125
\endthanks
\rightheadtext{II. Normal Cartan connections}
\endtopmatter

\document
In the sequel, we shall use the
notation and results from the preceding part of
this series of papers, see \cite{\v Cap, Slov\'ak, Sou\v cek}.
In particular, the citations like I.2.3
mean the corresponding items in that part.

In the first part we defined almost Hermitian symmetric structures as
`second order structures', see I.3.4. In this part we will first show that
any first order structure with the `right' structure group gives rise to
an almost Hermitian symmetric structure in this sense. Basically, the
construction is just the standard first prolongation of $G$--structures, see
\cite{Kobayashi} or \cite{Sternberg}. Due to the special situation there is
a canonical prolongation in this case which admits the structure of a
principal bundle with the structure group $B$.  Moreover, it turns out that
for the almost Hermitian structures in question, there exists a unique
normal Cartan connection. We shall prove this in section \nmb!{2}. Thus, the
calculus developed in part I, will yield natural operators in all these
cases.

Our approach is quite different to that of Ochiai whose vanishing torsion
assumption restricts in fact the considerations to the locally flat
structures in many cases, cf\. \cite{Baston, 91} or \cite{\v Cap, Slov\'ak,
95}. Another approach to the construction of the canonical Cartan
connections on certain auxiliary vector bundles, thus avoiding the
construction of the prolongation, can be found in \cite{Baston, 91}.

We shall discuss the applications to the
individual almost Hermitian symmetric structures in the next part
of this series.

\head \nmb0{1}. The  prolongation of first order
structures\endhead
\subheading{\nmb.{1.1}}
Let us recall the setting we work in from the first part of this
series of papers: We start from a connected semisimple real Lie group
$G$ whose Lie algebra $\goth g$ is equipped with a grading
$\goth g=\goth g_{-1}\oplus\goth g_0\oplus\goth g_1$. By $B$ we
denote the closed (parabolic) subgroup corresponding to the Lie
subalgebra $\goth b=\goth g_0\oplus\goth g_1$, further we have the
closed subgroup $B_0\subset B$  corresponding to $\goth g_0$ and the
closed normal subgroup $B_1$ of $B$ corresponding to $\goth g_1$.
Then it is known that
\roster
\item $\goth g_0$ is reductive with one--dimensional center
\item the map $\goth g_0\to\goth g\goth l(\goth g_{-1})$ induced by
       the adjoint representation is the inclusion of a subalgebra
       and an irreducible representation
\item the Killing form identifies $\goth g_1$ as a $\goth g_0$ module
       with the dual of $\goth g_{-1}$.
\item the restrictions of the exponential map to $\goth g_1$ and
       $\goth g_{-1}$ are diffeomorphisms onto the corresponding
       closed subgroups of $G$.
\item $B_0\cap B_1=\{e\}$
\item $B$ is the semidirect product of $B_0$ and $B_1$,
\endroster
see \cite{Ochiai, Sections 3 and 6}.

Examples of such Lie algebras and the corresponding structures can be
found in I.3.3. In particular recall that there are the classical
projective structures, which occur in this picture as the extremal
case of an almost Grassmannian structure. The projective structures
behave rather exceptionally and we will have to treat them separately.

\subheading{\nmb.{1.2}}
Our starting point is a first order $B_0$--structure on
a smooth manifold $M$ of dimension $m=\text{dim}(\goth g_{-1})$, so
assume that we have given a principal $B_0$ bundle $P_0\to M$
together with a {\it soldering form\/}
$\theta_{-1}\in \Om^1(P_0,\goth g_{-1})$ which is strictly
horizontal, i.e\. its kernel in each tangent space is precisely the
vertical tangent space, and $B_0$ equivariant, so
$(r^b)^*\theta_{-1}=\Ad(b^{-1})\o \theta_{-1}$. This is equivalent to
$P_0$ being a reduction of the (first order) frame bundle $P^1M$ of
$M$, cf\. I.3.6. Now consider the tangent bundle $TP_0$, the vertical
subbundle $VP_0$ and the quotient bundle $TP_0/VP_0$. The fundamental
vector field map gives a trivialization $VP_0\simeq P_0\x \goth g_0$,
while the soldering form induces a trivialization
$TP_0/VP_0\simeq P_0\x \goth g_{-1}$.

For a point $u\in P_0$ consider a linear isomorphism
$\ph:\goth g_{-1}\oplus \goth g_0\to T_uP_0$ which is compatible with
the two trivializations from above, i.e\. such that
$\ph(0,A)=\zeta_A(u)$ and $\theta_{-1}(u)(\ph(X,A))=X$.
Via $\ph$ the exterior derivative $d\theta_{-1}(u)$ gives rise to a
mapping $\goth g_{-1}\wedge\goth g_{-1}\to \goth g_{-1}$, defined by
$(X,Y)\mapsto d\theta_{-1}(u)(\ph(X,0),\ph(Y,0))$, and we view this
mapping as
$t_\ph\in\goth g_{-1}^*\wedge\goth g_{-1}^*\otimes\goth g_{-1}$,
and call it the {\it torsion\/} of $\ph$.
Now let $\bar\ph$ be another isomorphism compatible with the
trivializations. Then there is a linear map
$\ps:\goth g_{-1}\to \goth g_0$ such that
$\bar\ph (X,A)-\ph(X,A)=\zeta_{\ps(X)}(u)$. The difference between
the corresponding maps constructed using $d\theta_{-1}(u)$ can be
easily computed:
\proclaim{Lemma} In this situation we have:
$$d\theta_{-1}(u)(\bar\ph(X,0),\bar\ph(Y,0))-
d\theta_{-1}(u)(\ph(X,0),\ph(Y,0))=-[\ps(X),Y]+[\ps(Y),X].$$
\endproclaim
\demo{Proof}
Using bilinearity of $d\theta_{-1}(u)$ and the fact the
$\bar\ph(X,A)=\ph(X,A)+\zeta_{\ps(X)}(u)$ the difference can be
expressed as
$$d\theta_{-1}(u)(\zeta_{\ps(X)},\ph(Y,0))+
d\theta_{-1}(u)(\ph(X,0),\zeta_{\ps(Y)})+
d\theta_{-1}(u)(\zeta_{\ps(X)},\zeta_{\ps(Y)}).$$
Since $\theta_{-1}$ is horizontal and the Lie bracket of vertical
vector fields is vertical, the last term vanishes. On the other hand,
the infinitesimal version of the $B_0$--equivariancy of $\theta_{-1}$
is clearly $\Cal L_{\zeta_A}\theta_{-1}=-\ad(A)\o\theta_{-1}$, and
again by horizontality this reduces to
$i_{\zeta_A}d\theta_{-1}=-\ad(A)\o\theta_{-1}$. Applying this we see
that the first term from above reduces to
$-[\ps(X),\theta_{-1}(\ph(Y,0))]=-[\ps(X),Y]$ and similarly for the
second term.
\qed\enddemo

There is a canonical map $\partial$ from
$L(\goth g_{-1},\goth g_0)\simeq \goth g_{-1}^*\otimes\goth g_0$
to $\goth g_{-1}^*\wedge\goth g_{-1}^*\otimes \goth g_{-1}$,
the composition of the alternation in the first two
factors with the map induced by the inclusion
$\goth g_0\to\goth g_{-1}^*\otimes\goth g_{-1}$ obtained from
\nmb!{1.1}.\therosteritem{2}. Using this map, the lemma above just says
that
$$d\theta_{-1}(u)(\bar\ph(X,0),\bar\ph(Y,0))-
d\theta_{-1}(u)(\ph(X,0),\ph(Y,0))=-\partial(\ps)(X,Y).$$

Thus the above construction gives rise to a well defined function
$$P_0\to (\goth g_{-1}^*\wedge\goth g_{-1}^*\otimes\goth g_{-1})/
\partial(\goth g_{-1}^*\otimes\goth g_0),$$
called the {\it structure function\/} of the $B_0$--structure.

\subheading{\nmb.{1.3}}
The map $\partial:\goth g_{-1}^*\otimes\goth g_0\to
\goth g_{-1}^*\wedge\goth g_{-1}^*\otimes\goth g_{-1}$ from above is
the differential in the Spencer cohomology, the cohomology of the
abelian Lie algebra $\goth g_{-1}$ with values in the representation
$\goth g$. It is a crucial fact for the computation of this
cohomology that there is an adjoint
$\partial^*:\goth g_{-1}^*\wedge\goth g_{-1}^*\otimes\goth g_{-1}\to
\goth g_{-1}^*\otimes\goth g_0$, defined by
$(\partial^*\ph)(X)=\sum_i[Z^i,\ph(X_i,X)]$, where $\{X_i\}$ is a
basis of $\goth g_{-1}$ and $Z^i$ is the dual basis of $\goth g_1$,
see \nmb!{1.1}.\therosteritem{3}. It turns out that there is an
inner product on $\goth g$ such that $\partial^*$ is the adjoint
of $\partial$, see \cite{Ochiai, Proposition 4.2}. Thus the kernel
$\text{Ker}(\partial^*)$ is a complementary subspace to the image of
$\partial$.

Note that all spaces occurring in the above considerations are in fact
$\goth g_0$--modules. It is easy to verify that both $\partial$ and
$\partial^*$ are in fact homomorphisms of $\goth g_0$--modules. In
particular, this implies that $\text{Ker}(\partial^*)$ is even a
complementary $\goth g_0$--module to the image of $\partial$. This
will be crucial in the sequel.

\subheading{\nmb.{1.4}}
Now we define $P$ to be the set of all linear isomorphisms
$\ph:\goth g_{-1}\oplus\goth g_0\to T_uP_0$ as in \nmb!{1.2} such
that $\partial^*(t_\ph)=0$. It is easy to see that for each $u\in P_0$
such $\ph$ actually exist as follows: Take any $\ph$ satisfying the
conditions of \nmb!{1.2}. Then, as $\text{Ker}(\partial^*)$ is
complementary to $\text{Im}(\partial)$, there is a linear map
$\ps\in \goth g_{-1}^*\otimes\goth g_0$ such that
$\partial^*\bigl(d\theta_{-1}(u)(\ph(X,0),\ph(Y,0))+(\partial\ps)(X,Y)\bigr)=0$.
(In fact the image of $\ps$ under $\partial$ is uniquely determined.)
Then one immediately verifies that
$\bar\ph(X,A):=\ph(X,A)+\zeta_{\ps(X)}(u)$ satisfies the condition.

Next take an element $b\in B$. Viewing $b$ as an element of $G$ we
have the adjoint action $\Ad(b):\goth g\to \goth g$, and since
$\Ad(\text{exp}(Z))\cdot X=X+[Z,X]+1/2[Z,[Z,X]]+\dots$ (cf\.
I.3.8), we
see that $\goth g_1$ is stable under this adjoint action, so we get
an induced linear automorphism $\Ad(b)$ of the space
$\goth g/\goth g_1\simeq \goth g_{-1}\oplus \goth g_0$.

For $b\in B$ denote by $b_0$ the class of $b$ in $B/B_1\simeq B_0$.
Then for an element $\ph:\goth g_{-1}\oplus\goth g_0\to T_uP_0$ of
$P$ we define
$\ph\cdot b:\goth g_{-1}\oplus\goth g_0\to T_{u\cdot b_0}P_0$ by
$\ph\cdot b:=Tr^{b_0}\o\ph\o\Ad(b)$, where $r^{b_0}$ denotes the
principal right action of $b_0$ on $P_0$.

\proclaim{\nmb.{1.5}. Proposition}
This defines a free right action of $B$ on $P$. In each case except
the one of projective structures this action is also transitive on
each fiber of the obvious projection $P\to M$.
\endproclaim
\demo{Proof}
Let us first verify that $\ph\cdot b$ is again in $P$. So we have to
compute
$$
d\theta_{-1}(u\cdot b_0)((\ph\cdot b)(X,0),(\ph\cdot b)(Y,0)),
$$
for elements $X,Y\in \goth g_{-1}$. By $B_0$--equivariancy of
$\theta_{-1}$ this equals
$$
\Ad(b_0^{-1})(d\theta_{-1}(u)(\ph(\Ad(b)\cdot
(X,0)),\ph(\Ad(b)\cdot(Y,0)))).
$$
Now we may write $b=b_0b_1$ for some $b_1\in B_1$ and by
\nmb!{1.1}.\therosteritem{4} there is a $Z\in \goth g_{1}$ such that
$b_1=\text{exp}(Z)$. Using the formula for the adjoint action of an
exponential from above we see that
$$\Ad(b)\cdot(X,0)=\Ad(b_0)\cdot \Ad(\text{exp}(Z))(X,0)=
(\Ad(b_0)\cdot X,\Ad(b_0)\cdot [Z,X]),$$
and thus
$\ph(\Ad(b)\cdot(X,0))=\ph(\Ad(b_0)\cdot X,0)+\ze_{\Ad(b_0)\cdot [Z,X]}(u)$.
The same computation as in the proof of lemma \nmb!{1.2} then shows
that
$$\align
d\theta_{-1}(u)&(\ph(\Ad(b)\cdot (X,0)),\ph(\Ad(b)\cdot(Y,0)))=\\
&=d\theta_{-1}(u)(\ph(\Ad(b_0)\cdot X,0),\ph(\Ad(b_0)\cdot Y,0))+\\
&\newquad \Ad(b_0)\cdot ([[Z,X],Y]-[[Z,Y],X])\\
&=d\theta_{-1}(u)(\ph(\Ad(b_0)\cdot X,0),\ph(\Ad(b_0)\cdot Y,0)).
\endalign$$
This shows that $t_{\ph\cdot b}=b_0\cdot t_\ph$, so
$\partial^*(t_{\ph\cdot b})=b_0\cdot \partial^*(t_\ph)=0$, and hence
$\ph\cdot b\in P$.

Next, let us assume that $\ph\cdot b=\ph$ for some $\ph\in P$ and
$b\in B$. Then
obviously $b\in B_1$, since $B_0$ acts freely on $P_0$. So as before
we may write $b=\text{exp}(Z)$. But then $\ph\cdot b=\ph$ implies
that $[Z,X]=0$ for all $X\in\goth g_{-1}$, which implies $Z=0$ by
\nmb!{1.1}.\therosteritem{3}.

Finally, to prove transitivity of the action it suffices to show that
$B_1$ acts transitive on each fiber of $P\to P_0$, since $B_0$ acts
transitive on each fiber of $P_0\to M$. But for two maps
$\ph,\bar\ph$ in the same fiber we see from \nmb!{1.2} that
$\bar\ph(X,A)=\ph(X,A)+\zeta_{\ps(X)}(u)$ for some
$\ps\in\goth g_{-1}^*\otimes \goth g_0$, and lemma \nmb!{1.2} shows that
if both maps are in $P$ we must have $\partial(\ps)=0$. But now in all
cases except the projective one the corresponding Spencer cohomology
group $H^{1,1}(\goth g)$ is trivial, so there is a $Z\in \goth g_1$
such that $\ps=\ad_Z$, see \cite{Ochiai, Proposition 7.3}. Thus
$\bar\ph=\ph\cdot\text{exp}(Z)$.
\qed\enddemo

\subheading{\nmb.{1.6}. The soldering form}
>From now on we exclude the projective case which we will discuss
separately later. So $P\to M$ is a principal $B$--bundle, and the
proof of \nmb!{1.5} also shows that $p:P\to P_0$ is a principal
$B_1$--bundle. Now we define on $P$ a one--form $\theta$ with values
in $\goth g_{-1}\oplus \goth g_0$ as follows: For a point $\ph\in P$
consider a tangent vector $\xi\in T_\ph P$. Then $Tp\cdot \xi$ is a
tangent vector in $T_{p(\ph)}P_0$ and by definition $\ph$ is a linear
isomorphism from $\goth g_{-1}\oplus \goth g_0$ to this tangent
space, so we may define $\theta(\xi):=\ph^{-1}(Tp\cdot \xi)$. This
form is called the {\it soldering form\/} or {\it displacement
form\/} on $P$. The {\it torsion\/} $T$ of $\th$ is defined by the
structure equation
$$
d\th_{-1}=-[\th_0,\th_{-1}] + T.
$$

\proclaim{Lemma}
The one form $\theta$ has the following properties:
\roster
\item the component $\theta_{-1}$ is the pullback of the form from
       \nmb!{1.2}.
\item $\theta_0(\zeta_{Y+Z})=Y$ for all $Y\in \goth g_0$,
       $Z\in \goth g_1$.
\item $\theta$ is $B$--equivariant, i.e\.
       $(r^b)^*\theta=\Ad(b^{-1})\o\theta$, where $\Ad$ is the action
       from {\rm \nmb!{1.4}}
\item The torsion $T$ is horizontal over $M$ and can be viewed as a
function in $C^\infty(P,\goth g_{-1}^*\wedge\goth
g_{-1}^*\otimes \goth g_{-1})$. Moreover, $\partial^*\o T= 0$.
\endroster
In particular $(P,\theta)$ is a $B$--structure on $M$ in the sense of
{\rm I.3.4}.
\endproclaim
\demo{Proof}
\therosteritem{1} is clear since $\theta_{-1}(\ph(X,A))=X$. For
\therosteritem{2} note that $\zeta_Z$ lies in the kernel of $Tp$,
while $\zeta_Y$ is mapped by $Tp$ to the fundamental vector field on
$P_0$ corresponding to $Y$. Next, \therosteritem{3} follows
immediately from the definition of the $B$--action on $P$, and the
fact that $p\o r^b=r^{b_0}\o p$. Finally, $i_{\ze_X}$ with $X\in
\goth g_1$, applied to any
of the terms in the structure equation yields zero, while  for
$X\in \goth g_0$ we obtain
$$
i_{\ze_X}(d\th_{-1}+[\th_0,\th_{-1}]) = \L_{\ze_X}\th_{-1}
+[i_{\ze_X}\th_0,\th_{-1} ] = 0
$$
by the equivariancy of $\th_{-1}$. Now, we can define $T(X,Y)(u)$ by
evaluating the structure equation on arbitrary vectors $\xi$, $\et\in
T_uP$ such that $\th_{-1}(\xi)=X$ and $\th_{-1}(\et)= Y$. It remains
to prove $\partial^*\o T=0$ which can be done pointwise. So take
$\ph\in P$ and choose $\xi,\et\in T_\ph P$ so that $Tp.\xi= \ph(X,0)$
and $Tp.\et= \ph(Y,0)$. Then $\th_0(\xi)=\th_0(\et)=0$ by the
construction. Using (1), we see that
$$
d\th_{-1}(\xi,\et)+[\th_0,\th_{-1}](\xi,\et) = t_\ph(X,Y)+0
.\rlap{\hskip2cm\qed}$$
\enddemo

\subheading{\nmb.{1.7}}
Consider a principal connection $\gamma$ on $P_0$. Then at each point
$u\in P_0$ we get an isomorphism $\goth g_{-1}\oplus \goth g_0\to
T_uP_0$ as in \nmb!{1.2},  defined by the soldering form $\theta_{-1}$
and the connection form of $\ga$. Thus we have the torsion
$t_\ga:P_0\to \goth g_{-1}^*\wedge\goth g_{-1}^*\otimes\goth g_0$,
which is in fact the frame form of the usual torsion of $\ga$.

The connection $\ga$ is called {\it harmonic\/} if $\partial^*\o
t_\ga=0$.
(The name harmonic is due to the fact that the Spencer coboundary
operator $\partial$ is
trivial on the space in question so that our condition is equivalent
to harmonicity of the torsion.)

\proclaim{Proposition}
There is a $B_0$--equivariant section $\si:P_0\to P$, and the space of
all such sections is in bijective correspondence with the space of all
harmonic principal connections on $P_0$. Moreover, it is an affine
space modeled on $\Om^1(M)$, the space of one--forms on $M$.
\endproclaim
\demo{Proof} We have already shown in I.3.6 that a global
$B_0$-equivariant section $\si$ always exists, but now we shall
supply another simple (and more geometric) argument.

Note first that any principal connection $\ga$ on $P_0$ splits the
exact sequence
$$0\to VP_0\to TP_0\to TP_0/VP_0\to 0$$
and thus gives rise to
a linear isomorphism $\ph_u\:\goth g_{-1}\oplus \goth g_0 \to T_uP_0$,
which satisfies the conditions of \nmb!{1.2}, in
each point $u$.

Further, let us choose a $B_0$-module homomorphism $\ps$ which is a
right inverse of $\partial\: \goth
g_{-1}\otimes \goth g_0\to\text{Im}(\partial)\subset \goth
g_{-1}^*\wedge\goth g_{-1}^*\otimes \goth g_{-1}$.
Starting from a chosen principal connection $\ga$, let $f$ be the
$\text{Im}(\partial)$--part of the torsion $t_\ga$, and consider the
smooth map $u\mapsto \ph_u+\ze_{\ps(f(u))}(u)$. By the construction,
this has values in $P$, since by Lemma \nmb!{1.2} $t_{\ph+\ze_{\ps\o
f}} = t_\ph - \partial\o \ps\o f = t_\ph - f$. Due to the equivariancy
of $\ps$, this defines a $B_0$--equivariant section of $P\to P_0$. If
the original connection $\ga$ was harmonic, then $f=0$ and the mapping
$u\mapsto \ph_u$ itself is a $B_0$-equivariant section.

Any $B_0$-equivariant section $\si\:P_0\to P$ can clearly be
interpreted as a principal connection $\ga$ on $P_0$.
For each point $u\in P_0$ and $\xi\in T_uP_0$, we have
$$
(\si^*\th)(u)(\xi) = \th(\si(u))(T\si.\xi)= \si(u)^{-1}(\xi)\in \goth
g_{-1}\oplus \goth g_0
$$
and the $\goth g_0$--part of this expression is just the connection
form of the connection $\ga$. Applying $\si^*$ to the structure
equation from \nmb!{1.6}.(4) we obtain (using \nmb!{1.6}.(1))
$$
d\th_{-1} = -[\ga,\th_{-1}] + \si^* T,
$$
so that $\si^* T$ is the torsion of the principal connection
$\ga$. Thus $\ga$ is a harmonic connection.

Finally, if $\si$ and $\bar\si$ are two $B_0$-equivariant sections of
$P\to P_0$, then there is a unique smooth map $\ta\: P_0\to \goth g_1$
such that $\bar\si(u)= \si(u).\text{exp}\ta(u)$. Since the sections
are $B_0$-equivariant, we obtain $\ta(u.b_0)=\Ad(b_0\i).\ta(u)$, so
that $\ta$ is a frame form of a one-form on $M$.
\qed\enddemo

\subheading{\nmb.{1.8}}
The bundle $P$ can be viewed as a subbundle of the frame
bundle $P^1P_0$ of $P_0$. In fact, a point $\ph\in P$ is by
definition an isomorphism
$\goth g_{-1}\oplus \goth g_0\to T_{p(\ph)}P_0$. Moreover, taking
into account that $P_0$ is a reduction of $P^1M$ to the group $B_0$,
we can view $P$ as a reduction of $P^1(P^1M)$ to the group $B$. In
fact, it can be shown that this reduction has values in the second
order frame bundle $P^2M$ of $M$, if and only if the torsion of
$\theta$ vanishes, cf\. \cite{Slov\'ak, 94}, but we will not pursue
this point of view.

\subheading{\nmb.{1.9}. The projective case}
In this case the underlying first order structure is the whole
$P^1M$, so it carries no information. Thus to get a $B$--structure in
the sense of I.3.4 with harmonic torsion, one has to choose a
reduction of the second order frame bundle $P^2M$ to the appropriate
group $B$. (Note that in this case
$\partial:\goth g_{-1}^*\otimes \goth g_0\to
\goth g_{-1}^*\wedge\goth g_{-1}^*\otimes \goth g_{-1}$ is
surjective, so the harmonic connections are exactly the torsion free
ones.)

\head\nmb0{2}. Canonical Cartan connections\endhead
Our next task is to prove that in all but the very low dimensional
cases, on all prolongations as constructed in the previous  section,
there is a canonical Cartan connection. Basically, this is a
consequence of the fact that in these cases the next prolongation is
trivial, so its soldering form is a Cartan connection.

\subheading{\nmb.{2.1}}
Assume we have constructed the $B$--bundle $P\to M$ with the soldering
form $\theta=\th_{-1}\oplus\th_0$ for a $B_0$--bundle $P_0\to M$ as
above. As we have seen in \nmb!{1.8} this is in fact a $B_1$ structure
on $P_0$, so we can try to apply the same construction as above to
this structure using the additional information we have in this case.

The starting point is to consider for $\ph\in P$ linear isomorphisms
$$\Ph:\goth g_{-1}\oplus \goth g_0\oplus \goth g_1\to T_\ph P$$
such that $\th(\Ph(X,A,Z))=(X,A)$ and $\Ph(0,A,Z)=\ze_{A+Z}(\ph)$ (here we
use the finer structure and do not only  fix $\Ph(0,0,Z)$). Having given
such a $\Ph$ we have to consider its torsion
$$\gather
t_\Ph\in (\goth g_{-1}\oplus \goth g_0)^*\wedge(\goth g_{-1}\oplus \goth g_0)^*
\otimes(\goth g_{-1}\oplus \goth g_0)\\
((X,A),(Y,B))\mapsto d\th(\ph)(\Ph(X,A,0),\Ph(Y,B,0)).
\endgather$$
In fact, several parts of this mapping are independent of $\Ph$. For
later use we prove a slightly more general result than we need here:
\proclaim{Lemma}
Let $pr:\goth g\to \goth g_{-1}\oplus \goth g_0$ be the obvious
projection. Then for all $(X,A,Z)$, $(Y,B,W)\in \goth g$ we have
$$
d\th(\Ph(X,A,Z),\Ph(Y,B,W))=d\th(\Ph(X,0,0),\Ph(Y,0,0))-pr([X+A+Z,Y+B+W]).
$$
\endproclaim
\demo{Proof}
The infinitesimal version of the equivariancy of $\th$ gives
$$i_{\ze_{A+Z}}d\th=-\ad(A+Z)\o \theta,$$
where $\ad$  is the
composition of $pr$ with the adjoint action on $\goth g$. Now the
result follows easily using bilinearity of $d\th$ and the fact that
$\Ph(X,A,Z)=\Ph(X,0,0)+\ze_{A+Z}(\ph)$.
\qed\enddemo
Consequently, the torsion of $\Ph$ is determined by its component
in $\goth g_{-1}^*\wedge \goth g_{-1}^*\otimes \goth g_0$.

\subheading{\nmb.{2.2}}
The next step is to compute the change of the torsion if one replaces
$\Ph$ by another isomorphism satisfying the above conditions. As in
the proof of \nmb!{1.2} one verifies that in fact the change lies in
the image of $\goth g_{-1}^*\otimes \goth g_1$ under
the composition of the alternation with the map induced by the
inclusion
$$\goth g_{-1}^*\otimes \goth g_1\to \goth g_{-1}^*\otimes
\goth g_{-1}^*\otimes \goth g_0\subset (\goth
g_{-1}\oplus \goth g_0)^*\otimes (\goth
g_{-1}\oplus \goth g_0)^*\otimes (\goth g_{-1}\oplus \goth g_0).$$
To get a well defined structure function as in section \nmb!{1} we
have to factor the latter space by the image of
$\partial:\goth g_{-1}^*\otimes \goth g_1\to \goth g_{-1}^*\wedge \goth
g_{-1}^*\otimes \goth g_0$. As before this is the differential in the
Spencer cohomology and it has an adjoint $\partial^*$ defined by the
same formula as in \nmb!{1.3}.

\proclaim{\nmb.{2.3}. Theorem}
In all cases but the one of $\goth g=\goth s\goth l(2)$, for each
$\ph\in P$ there is a unique linear isomorphism $\Ph$ as in \nmb!{2.1}
such that $\partial^*\o t_\Ph=0$. The inverses of these can be viewed
as a smooth one form $\om\in \Om^1(P,\goth g)$ with the following properties:
\roster
\item $\om(\ze_X) = X$ for all $X\in \goth b$
\item $(r^b)^*\om=\text{Ad}(b^{-1})\o \om$ for all $b\in B$
\endroster
\endproclaim
\demo{Proof}
First, since the kernel of $\partial^*$ is a complement to the image
of $\partial$, we can construct such a $\Ph$ in a point $\ph$ like in
\nmb!{1.4}. Moreover, it is clear that the set of all such $\Ph$ is
parameterized by the kernel of $\partial:\goth g_{-1}^*\otimes \goth
g_1\to \goth g_{-1}^*\wedge \goth g_{-1}^*\otimes\goth g_0$. This
coincides with the Spencer cohomology group $H^{2,1}(\goth g)$ which
is trivial for all cases in question, see \cite{Ochiai, Proposition
7.1}, so $\Ph$ is unique.

Let us verify the properties of $\om$. For $A\in\goth g_0$ and
$Z\in\goth g_1$ we have
$\Ph(0,A,Z)=\zeta_{A+Z}(\ph)$, so $\om$ reproduces the generators of
fundamental vector fields. Finally, we have to verify the equivariancy of
$\om$. Put $\Ph=\om(\ph)^{-1}:\goth g\to T_\ph P$ and consider
$\Ph\cdot b:=Tr^b\o \Ph\o \Ad(b):\goth g\to T_{\ph\cdot b}P$ for $b\in
B$. If we
verify that $\Ph\cdot b$ satisfies the conditions of \nmb!{2.1} and
that $\partial^*\o t_{\Ph\cdot b}=0$, then the uniqueness proved above
concludes the proof.

For $A\in \goth g_0$ and $Z\in\goth g_1$ we have
$$(\Ph\cdot
b)(0,A,Z)=Tr^b\zeta_{\Ad(b)\cdot (0,A,Z)}(\ph)=\zeta_{A+Z}(\ph\cdot
b).$$

Further, $\theta(\ph\cdot b)((\Ph\cdot
b)(X,A,Z))=\Ad(b^{-1})\th(\ph)(\Ph(\Ad(b)\cdot
(X,A,Z)))=(X,A)$, since $\Ad(b)\cdot
(X,A)$ is by definition just the first two components of $\Ad(b)\cdot
(X,A,Z)$.

It remains to check the condition on the torsion. For $b\in B$,
we write $b=b_0\text{exp}(W)$ (see \nmb!{1.5}). Using the
equivariancy of $\theta$ and lemma \nmb!{2.1} we compute:
$$\align
d\th&(\ph\cdot b)(Tr^b\Ph(\Ad(b)\cdot X),Tr^b\Ph(\Ad(b)\cdot Y))=\\
&=\Ad(b^{-1})(d\th(\ph)(\Ph(\Ad(b)\cdot X),\Ph(\Ad(b)\cdot Y)))\\
&=\Ad(b^{-1})(d\th(\ph)(\Ph(\Ad(b_0)\cdot X),\Ph(\Ad(b_0)\cdot Y)))+
pr([\Ad(b)\cdot X,\Ad(b)\cdot Y])
\endalign$$
The second term in this expression vanishes since $\Ad(b)$ is a Lie
algebra homomorphism, so for the $\goth g_0$--component we get
$$\multline
\Ad(b_0^{-1})(d\th_0(\ph)(\Ph(\Ad(b_0)\cdot X),\Ph(\Ad(b_0)\cdot Y)))+\\
+\Ad(b_0^{-1})([W,d\th_{-1}(\ph)(\Ph(\Ad(b_0)\cdot X),\Ph(\Ad(b_0)\cdot Y))]).
\endmultline$$
The first term lies in the kernel of $\partial^*$ since this is a
$B_0$--submodule. The second one lies in this kernel since by
definition of $\partial^*$ we have
$\partial^*\o\ad(W)=\ad(W)\o\partial^*$ (cf\. \nmb!{1.3}).
\qed\enddemo

Since the restriction of the one form $\om$ to any tangent space
$T_\ph P$ is an isomorphism, $\om$ is a Cartan connection on $P$, see
the definition in I.2.1. Moreover, the first condition put on $\Ph$ in
\nmb!{2.1} implies that the $(\goth g_{-1}\oplus\goth g_0)$--part of
$\om$ coincides with $\th$. Thus, $\om$ is an admissible Cartan
connection in the sense of I.3.9.

Let us remark that another approach to the construction of canonical
prolongations equipped with canonical Cartan connections can be found
in \cite{Alekseevsky, Michor, 93}.

\subheading{\nmb.{2.4}}
Let us return to the point of view of $G$--strucutres and compute the
structure function of the last prolongation. Clearly this is induced
by
$$
(X,A,Z),(Y,B,W)\mapsto d\om(\om^{-1}(X,A,Z),\om^{-1}(Y,B,W)).
$$
Using the $B$--equivariancy of $\om$ one proves precisely as in lemma
\nmb!{2.1} that
$$
\multline d\om(\om^{-1}(X,A,Z),\om^{-1}(Y,B,W))=\\=d\om(\om^{-1}(X,0,0),\om^{-1}(Y,0,0))+
[X+A+Z, Y+B+W].\endmultline
$$

By definition the first term of the right hand side is just the
curvature $\ka(X,Y)$, see I.2.1. Thus, if this curvature vanishes
then the structure function is constant (and equal to the Lie
bracket,  viewed as an element of
$(\frak g^*\wedge \frak g^*)\otimes \frak g$), independent of the
manifold under consideration.

In this situation, taking into account that the
components $\om_{-1}$ and $\om_0$ coincide with the respective
components of $\th$, we see from \nmb!{2.1} that also the next
``lower'' structure function is constant and independent of the
manifold. Similarly, one shows that the same is true for the first
structure function constructed in \nmb!{1.2} .

In the flat case $M=G/B$ the canononical Cartan connection is just
the Maurer Cartan form, and the Maurer Cartan equation means just
that  $\ka=0$ in this case. Thus we see that a $B_0$--structure
$P_0\to M$ has the structure functions of all prolongations constant
and equal to those of the flat model if and only if the curvature of
the canonical Cartan connection vanishes. From \cite{Sternberg, p. 339}
we conclude:

\proclaim{Proposition}
$P_0\to M$ is locally isomorphic to the flat model if and only if the
canonical Cartan connection has zero curvature.
\endproclaim

\subhead\nmb.{2.5} \endsubhead
Using the properties of Cartan connections derived in \cite{\v Cap,
Slov\'ak, Sou\v cek} it is quite easy to compute explicitely the
obstructions against flatness of the canonical Cartan connection in
terms of any of the underlying linear connections. The main step is
to exploit the known results on the second cohomology
$H^2(\frak g_{-1},\frak g)$ to determine, which parts of the
curvature are true obstructions, and which vanish automatically.
Then one can explicitly compute the relevant components. This is
worked out in \cite{\v Cap, Slov\'ak}

\head\nmb0{3}. Explicit formulae for the canonical Cartan connections\endhead

Let us consider a $B$-structure $P\to P_0\to M$, its soldering form
$\th$ with a harmonic torsion $T$, and the canonically defined Cartan
connection $\om$ on $P$,
as constructed in section \nmb!{2}. Note that the canonical Cartan
connection is characterized by the fact the the component $\ka_0$ of
its curvature is in the kernel of $\partial^*$.

For each global
$B_0$--equivariant section $\si\:P_0\to P$ there is the principal
connection $\si^*\th_0$ on $P_0$, the induced
admissible Cartan connection $\tilde \ga$ on $P$,
and the difference between the
canonical Cartan connection and the latter one is described by the so
called deformation tensor $\Ga$, see I.3.9. In this section, we shall
compute explicitly the deformation of a chosen induced
admissible Cartan connection which leads to the canonical one. It
turns out, that for each of the structures in question, there is a
universal formula for $\Ga$ in terms of the curvature tensor of the
chosen underlying connection $\ga$.

Since the computations are quite elementary and in fact an explicit
use of the general result from section \nmb!{2} does not spare much
work, we prefer to recover completely also
the existence and uniqueness of the canonical Cartan connection in
this way. Thus a part of the next considerations will be redundant,
but on the other hand, this will also provide the link to the
traditional concept of the normal Cartan connection, see e.g\.
\cite{Kobayashi, 72}.

\subheading{\nmb.{3.1}. The trace of the curvature}
Let $\om$ be an admissible Cartan connection on a $B$-structure
$P\to M$, i.e\. $\om=\th_{-1}\oplus\th_0\oplus \om_1$.
Let us recall the definition of the trace of the
$\goth g_0$-component $\ka_0$ of the curvature function $\ka$ of
$\om$.  We can view the values of $\ka_{0}$ as elements
in $\goth g_{-1}^*\otimes\goth g_{-1}^*\otimes \goth
g_{-1}^*\otimes \goth g_{-1}=
\goth g_{-1}^*\otimes\goth g_{-1}^*\otimes \goth
g_{-1}^*\otimes \goth g_{1}^*$.

There are three
possible evaluations in the target space. The evaluation
over the last two entries is just the trace in $\goth g_{0}$,
the other two possibilities coincide up to a sign. By definition,
the \idx{\it trace} $\operatorname{Tr}\ka_0$ of the curvature
function $\ka_0$ is the evaluation over the first and the last
entry.

\proclaim{Lemma} For all $X\in \goth g_{-1}$ we have
$(\partial^{*}\ka_0)(X)= (\operatorname{Tr}\ka_0)(X,\_)\in \goth
g_1$. In particular, $\operatorname{Tr}\ka_0=0$ if and only
if  $\partial^* \ka_0=0$.
\endproclaim
\demo{Proof} By the definition above, $(\operatorname{Tr}\ka_0)(X,Y)
= \sum_i \ka_0(e_i,X)(Y)(e^i)$, where $e_i$ is a basis in $\goth
g_{-1}$ while $e^i$ is its dual basis in $\goth g_1$. If we take $Y$
as a free argument, we obtain $\operatorname{Tr}\ka_0(X,\_)\in
\goth g_1$,
$\operatorname{Tr}\ka_0(X,\_)=\sum_i[e^i,\ka_0(e_i,X)]$. But the
latter is exactly the formula for $(\partial^*\ka_0)(X)$, see
\nmb!{1.3}. \qed\enddemo

\subheading{\nmb.{3.2}. Definition}
A \idx{\it normal Cartan connection} $\om\in\Om^1(P,\goth g)$
is an admissible connection with the curvature satisfying
$\operatorname{Tr}\ka_0=0$.

\proclaim{\nmb.{3.3}. Lemma} Let $P\to M$ be a $B$ structure
with harmonic torsion, $P_0$ be the underlying first order
structure. Then for each admissible Cartan
connection $\om$ on $P$, $\partial^*\ka_0$ is constant on the
fibers of $P\to P_0$.
\endproclaim
\demo{Proof}
By the formula I.3.8.(4), for each  section $\si $ of $P\to P_0$
and $u\in P$ we have
$$\ka_0(u)(X,Y) = \ka_0(\si(p(u)))(X,Y)-[\ta(u),\ka_{-1}(\si(p(u)))(X,Y)],
$$
where $\tau$ is the mapping introduced in the proof of I.3.7.

Further, $\partial^*\ka_{-1}=0$ since the torsion $\ka_{-1}$ is
harmonic, and we obtain
$$\align
\partial^*\ka_0(u)(X) &= [e^i, \ka_0(\si(p(u)))(e_i,X)]-
[e^i,[\ta(u),\ka_{-1}(\si(p(u)))(e_i,X)]]\\
&= \partial^*\ka_0(\si(p(u)))(X) -
[\ta(u),[e^i, \ka_{-1}(\si(p(u)))(e_i,X)]]\\
&= \partial^*\ka_0(\si(p(u)))(X) - [\ta(u), \partial^*\ka_{-1}(\si(p(u)))(X)]
\\&= \partial^*\ka_0(\si(p(u)))(X).\qed
\endalign$$
\enddemo

We shall also need another technical lemma. In view of lemma
\nmb!{3.1}, it is a direct consequence of the uniqueness result
from the previous section, the elementary argument used here
is an easy application of the Bianchi identity for general Cartan
connections.

\proclaim{\nmb.{3.4}. Lemma} Let $\om$ and $\bar \om$ be two
normal Cartan connections on a $B$-structure $P$, $\ka$ and $\bar\ka$
be their curvatures. Then the trace
$\operatorname{Tr}_{\goth g_0}(\bar\ka_0-\ka_0)$ within $\goth
g_0$ vanishes.
\endproclaim
\demo{Proof}
Let us write $\de= \bar\ka_0-\ka_0$, and let $e_i$ and $e^i$ be the dual
bases in $\goth g_{\pm1}$. According to the Bianchi identity
(proved in I.2.4) we have for $X,Z\in\goth g_{-1}$
$$\align
[\de(X,Z)&,e_i] = [\de(X,e_i), Z] +[\de(e_i,Z),X] +
\\
&\nabla^{\om}_{Z}\ka_{-1}(X,e_i) +
\nabla^{\om}_{X}\ka_{-1}(e_i,Z)+\nabla^{\om}_{e_i}\ka_{-1}(Z,X) +
\\
& \ka_{-1}(\ka_{-1}(X,e_i),Z) +
\ka_{-1}(\ka_{-1}(e_i,Z),X)+\ka_{-1}(\ka_{-1}(Z,X),e_i)-
\\
&\nabla^{\bar\om}_{Z}\bar\ka_{-1}(X,e_i) -
\nabla^{\bar\om}_{X}\bar\ka_{-1}(e_i,Z)-\nabla^{\bar\om}_{e_i}
\bar\ka_{-1}(Z,X) -
\\
&\bar\ka_{-1}(\bar\ka_{-1}(X,e_i),Z) -
\bar\ka_{-1}(\bar\ka_{-1}(e_i,Z),X)-\bar\ka_{-1}(\bar\ka_{-1}
(Z,X),e_i).
\endalign$$
Since $\bar\ka_{-1}=\ka_{-1}$ and the torsion $\ka_{-1}$ is
constant on the fibers of $P\to P_0$, all lines except the first
one vanish, see I.3.8.(4), I.3.10.(3) and the definition of
$\nabla^{\om}$ in I.2.3.

\noindent Now, $\operatorname{Tr}_{\goth g_0}(\de)(X,Z) =
\sum_i[\de(X,Z),e_i](e^i)$ while $(\operatorname{Tr}\de)(X,Z) =
[\de(e_i,X), Z](e^i)=0$. Thus the above computation shows that
the traces inside of $\goth g_0$ coincide as required.
\qed\enddemo

\subheading{\nmb.{3.5}. Remark} If the torsion of a
$B$-structure $P$ vanishes, then all the admissible
Cartan connections have vanishing $\goth g_{-1}$-part of the
curvature. Then the Bianchi identity implies directly that
$\partial \ka_{0}$ vanishes for all admissible Cartan
connections. Thus,
in the language of the Hodge theory for the corresponding
cohomologies, this means just that the normal Cartan connections
are exactly those admissible Cartan connections
for which $\ka_0$ is harmonic.
As discussed in \nmb!{1.8}, if there is a torsion
free connection on a reduction $P_0$ of
$P^1M$ to the structure group $B_0$, then there is the canonical
$B$-structure $P$ over $P_0$ with vanishing torsion and a normal
Cartan connection on $P$ is then an admissible Cartan connection with
a harmonic $\goth g_0$-part of the curvature. This is
the point of view adopted in \cite{Ochiai, 70} where the
torsion--free case is discussed. However this
cannot yield a canonical Cartan connection in the cases of non vanishing
torsion in view of the results of the previous section.

\subhead \nmb.{3.6}. The conformal case\endsubhead
Since there is always a torsion--free linear connection on each
Riemannian manifold, the canonical prolongation $P$ of a first
order conformal Riemannian structure $P_0\to M$ is always a
reduction of $P^2M$ and so we reproduce the classical
construction in this case, cf\. \cite{Kobayashi, 72}.
We already deduced in I.6.3 the existence and uniqueness of the normal
Cartan connections, and the corresponding  explicit deformation
tensors $\Ga$, on all manifolds $M$ with conformal Riemannian
structures, $\operatorname{dim}M\ge3$. Let us recall the final
formula: Starting with a torsion--free connection $\ga$ on $P_0$
with curvature tensor $R^i_{jkl}$, Ricci tensor $R_{ij}$
and scalar curvature $R$, the necessary deformation
tensor $\Ga\in C^\infty(P_0, \goth g_{-1}^*\otimes\goth g_{-1}^*)$
is given by
$$
\Ga_{ij}=\tfrac{-1}{m-2}\bigl(R_{ij} -
\tfrac{\de_{ij}}{2(m-1)}R\bigr).
$$

\subhead \nmb.{3.7}. The almost Grassmannian case \endsubhead
Now we shall construct the normal Cartan connections on manifolds with
almost Hermitian symmetric structures corresponding to the
algebras $\goth g=\goth s\goth l(p+q, \Bbb R)$.
The description of the algebra
$\goth g=\text{Mat}_{q,p}(\Bbb R)\oplus (\goth s\goth
l(p,\Bbb R)\oplus \goth s\goth l(q,\Bbb R)\oplus \Bbb
R)\oplus \text{Mat}_{p,q}(\Bbb R)$ yields easily the formulas for
the bracket. Let us use the generators $e^\al_\be$ of the vector
spaces of matrices, the matrices
with all entries zero except a 1 in the $\be$--th line and the
$\al$--th column. We
shall use the letters $a,b,c,\dots$ for the indices between 1 and
$p$, the letters $i,j,k,\dots$ will indicate
indices running between 1 and $q$. For example, $e^a_i$ means one
of the generators in $\text{Mat}_{q,p}(\Bbb R)$. Using the fact that
the Killing form of $\goth g=\goth s\goth l(p+q, \Bbb R)$ is a scalar
multiple of the trace form one easily see that the bases $\{e^i_a\}$
and $\{e^a_i\}$ are also dual with respect to the Killing form, up to
a fixed scalar multiple, and this suffices for our purposes. Then we have
$$
[e^a_i, e^j_b] = \de^a_be^j_i - \de ^j_ie^a_b,~
[e^k_a, e^b_c] = -\de^b_ae^k_c,~[e^k_a,e^j_l] =
\de^k_le^j_a .
$$

Let us fix the sizes $p$ and $q$, and
consider an almost Grassmannian structure $P\to M$ with a harmonic
torsion. Let $P_0\to M$ be the underlying first order structure with the
distinguished class of the harmonic connections.

The deformation tensor $\Ga$ is expressed through functions
$\Ga_{{}^b_j{}^a_i}$ defined by $\Ga(e^a_i)=\Ga_{{}^b_j{}^a_i}e^j_b$. The
possible deformations $\de\ka_0$ of the curvature are described
in I.3.10.(4) The trace of the curvature is obtained through evaluation
in $\goth g_{-1}^*\otimes \goth g_{-1}^*\otimes \goth g_{-1}^*\otimes
\goth g_{-1}$ over the first and the fourth entry, however
according to Lemma \nmb!{3.1}, we
can compute $\partial^*(\ka_0)$ instead.
Let us first evaluate $[e^i_b, [\Ga.e^b_i,Y]-[\Ga.Y,e^b_i]]$ on the
generators.
$$\align
[e^i_b,[\Ga.e^a_k,e^b_i] - [\Ga.e^b_i, e^a_k]] &=
[e^i_b,[\Ga_{{}^d_s{}^a_k}e^s_d,e^b_i]-
[\Ga_{{}^d_s{}^b_i}e^s_d, e^a_k]]\\
&= [e^i_b,-\Ga_{{}^b_s{}^a_k}e^s_i+\Ga_{{}^d_i{}^a_k}e^b_d+
\Ga_{{}^a_s{}^b_i}e^s_k -\Ga_{{}^d_k{}^b_i}e^a_{d}]\\
&=(-\de^i_i\Ga_{{}^b_s{}^a_k}+\de^i_k\Ga_{{}^a_s{}^b_i})e^s_b +
(-\de^b_b\Ga_{{}^d_i{}^a_k}+\de^b_a\Ga_{{}^d_k{}^b_i})e^i_d
\endalign$$
The application of the formula for $\partial^*$ from \nmb!{1.3}
yields
$$\align
\partial^*(\de\ka_0)(e^a_k) &= \tsum_{s=1}^q\tsum_{b=1}^p
(-q\Ga_{{}^b_s{}^a_k}+\Ga_{{}^a_s{}^b_k})e^s_b +
\tsum_{d=1}^p\tsum_{i=1}^q
(-p\Ga_{{}^d_i{}^a_k}+\Ga_{{}^d_k{}^a_i})e^i_d
\tag1\\
&=\tsum_{l=1}^q\tsum_{c=1}^p
(-q\Ga_{{}^c_l{}^a_k} - p\Ga_{{}^c_l{}^a_k} + \Ga_{{}^a_l{}^c_k} +
\Ga_{{}^c_k{}^a_l})e^l_c.
\endalign$$
According to \nmb!{3.1}, the trace of $\ka_0$ evaluated on
the base elements $e^a_k$, $e^c_l$ is exactly the expression
inside the brackets in the last sum.

For each admissible Cartan connection $\om$ on $P$ there are the two
parts $\ka_{0,1}$, $\ka_{0,2}$ of $\ka_0$, corresponding to the
decomposition of $\goth g_0$ into two components. They are
given by functions $K^{a}_{d{}^b_k{}^c_l}$,
$K^{i}_{j{}^b_k{}^c_l}$, one set for each of the two blocks in
the matrices in $\goth g_0$.
>From the second line of the above computation, we can read the
formulae for the deformation of these functions achieved by the
chosen deformation tensor $\Ga$
$$
\de K^{l}_{s{}^b_i{}^a_k}=
\Ga_{{}^b_s{}^a_k}\de^l_i-\Ga_{{}^a_s{}^b_i}\de
^l_k,\newquad
\de K^{d}_{c{}^b_i{}^a_k}= \Ga_{{}^d_k{}^b_i}\de^a_c -
\Ga_{{}^d_i{}^a_k}\de^b_c
.$$
Consequently, the deformation of
the traces $\text{Tr}_{\goth g_0}(\de(\ka_{0,1}))$,
$\text{Tr}_{\goth g_0}(\de(\ka_{0,2}))$ of these two
components within $\goth g_0$ are
$\mp (\Ga_{{}^a_k{}^b_i}-\Ga_{{}^b_i{}^a_k})$.

Now, given a connection $\ga$ in the distinguished class on
$P_0$, we
shall compute the deformation tensor $\Ga$ which deforms the induced
admissible Cartan connection $\tilde \ga$ into a normal Cartan
connection $\om$ with curvature $\bar\ka = \ka -\de\ka$.
Let $\ka$ be the curvature function of $\tilde
\ga$, and write $\de\ka$ for its
change achieved by the choice of $\Ga$. We have
$$\align
\text{Tr}_{\goth g_0}(\de(\ka_{0,2}))_{{}^c_l{}^a_k}&=
\Ga_{{}^c_l{}^a_k}-\Ga_{{}^a_k{}^c_l}\\
(p+q)\text{Tr}(\de(\ka_0))_{{}^c_l{}^a_k}&=
-(p+q)^2\Ga_{{}^c_l{}^a_k} + 2(p+q)\Ga_{{}^a_l{}^c_k} -
(p+q)\text{Tr}_{\goth g_0}(\de(\ka_{0,2}))_{{}^a_l{}^c_k}\\
2\text{Tr}(\de(\ka_0))_{{}^a_l{}^c_k}&=
-2(p+q)\Ga_{{}^a_l{}^c_k} + 4\Ga_{{}^c_l{}^a_k} -
2\text{Tr}_{\goth g_0}(\de(\ka_{0,2}))_{{}^c_l{}^a_k}
\endalign$$
where the aim of our manipulation is to get rid of the
interchanging indices in the formula for the trace of $\ka_0$.

Let $\si\:P_0\to P$ be the section corresponding to the
connection $\ga$. The curvature $R=R_1+R_2\:P_0\to \goth g_{-1}^*\wedge\goth
g_{-1}^*\otimes \goth g_0$ of $\ga$ is $\si$-related to
$\ka_0=\ka_{0,1}+\ka_{0,2}$. In particular, on the image
$\si(P_0)\subset P$, we can achieve the vanishing of the trace
of $\bar\ka_0$ by the following choice of the deformation
$$
\Ga_{{}^c_l{}^a_k} = \tfrac{-1}{4-(q+p)^2}
\bigl((p+q)\text{Tr}(R)_{{}^c_l{}^a_k} +
2\text{Tr}(R)_{{}^a_l{}^c_k} + (p+q)\text{Tr}_{\goth g_0}(R_2)_{{}^a_l{}^c_k}
+ 2\text{Tr}_{\goth g_0}(R_2)_{{}^c_l{}^a_k}\bigr).
\tag2$$

Since the torsion of the $B$-structure $P$ is harmonic, this
choice of $\Ga$ leads to a normal Cartan connection according to
Lemma \nmb!{3.3}.

The results of the previous section assure that there is a unique
normal Cartan connection on $P$, however it is easy to verify
this directly. Indeed, it is equivalent to
prove, that if $\om$ and $\bar \om$ are two normal
Cartan connections on $P$, then the (uniquely defined) deformation
tensor $\Ga$ is identically zero.
In fact, we have computed above a tensor $\Ga$ deforming a given
$\om$ in such a way, that on the image of a section $\si\:P_0\to
P$ the achieved deformation of the trace of $\goth g_0$-part of
the curvature of $\ga$ reaches a value prescribed in advance. But
Lemma \nmb!{3.4} states that the traces of $\ka_0$ and
$\bar\ka_0$ inside of $\goth g_0$ coincide. In view of our
computation this means, that
the deformation tensor $\Ga$ satisfies
$\Ga_{{}^a_k{}^b_i}=\Ga_{{}^b_i{}^a_k}$ and so
for any two normal Cartan connections $\om$, $\bar\om$, the
corresponding deformation tensor $\Ga$ is symmetric.
Further, the achieved deformation of the trace of $\ka_0$ by means of
$\Ga$ has to vanish too and
since we can use the equality $\Ga_{{}^a_k{}^b_i}= \Ga_{{}^b_i{}^a_k}$
we obtain  $2\Ga_{{}^c_k{}^a_l}= (p+q)\Ga_{{}^c_l{}^a_k}$. Applying the
latter equality twice, we get
$$
2(p+q)\Ga_{{}^c_l{}^a_k} = 4\Ga_{{}^c_k{}^a_l} =
(p+q)(p+q)\Ga_{{}^c_k{}^a_l}
.$$
Thus, if $q\ge p\ge 1$, $q+p\ge 3$ then $\Ga_{{}^c_k{}^a_l}=0$
for all $c,a,k,l$ and so there is at most one normal Cartan
connection $\om$ on $P$.

Thus, we can formulate the final result of our computations.

\proclaim{\nmb.{3.8}. Theorem} Let $P\to M$ be a real
almost Grassmannian structure
with a harmonic torsion, on a smooth manifold $M$ and assume $q\ge
p\ge1$, $q+p\ge 3$.
Then there is a uniquely defined normal Cartan connection $\om$
on $P$ and for each linear harmonic connection $\ga$ on the
underlying first order structure $P_0$ with curvature $R=R_1+R_2$,
$\om = \tilde\ga-\Ga\o\th_{-1}$, where the corresponding deformation tensor
$\Ga$ is given by the formula {\rm \nmb!{3.7}.(2)}.
\endproclaim

\proclaim{\nmb.{3.9}. Corollary} Let $P\to M$ be a projective
structure on a smooth manifold $M$, $\operatorname{dim}(M)=q>1$.
Then there is a uniquely defined normal Cartan connection $\om$
on $P$ and for each linear torsion--free connection $\ga$ from the
underlying class on the first order structure $P_0$ with curvature
$R=(R^i_{jkl})$, we obtain
$\om = \tilde\ga-\Ga\o\th_{-1}$, where the corresponding deformation tensor
$\Ga$ is given by
$$
\Ga_{jk} = \tfrac{1}{(q-1)}(R^l_{jlk}+ R^l_{ljk}).
$$
\endproclaim

\subheading{\nmb.{3.10}. The almost Lagrangian case} We have to
deal with $\goth g = \goth g_{-1}\oplus\goth g_0\oplus\goth g_1$
where $\goth g_{-1}= S^2\Bbb R^m$, $\goth g_1=S^{2}\Bbb R^{m*}$,
$\goth g_0=\goth g\goth l(m,\Bbb R)$, cf\. I.3.3. Let us
fix the base $e_k\odot e_l$ consisting of symmetric
matrices with entries
$a^i_{j}=\frac12(\de^i_k\de_{jl}+\de^i_l\de_{jk})$.
Let $e^s\odot e^t$ be the dual base of $\goth g_1$ and let $e^i_j$
be the usual base of $\goth g_0$. The commutators of
the base elements are
$$\align
[e^s\odot e^t, e_k\odot e_l] &= - \frac14(\de^s_ke^t_l +
\de^s_le^t_k + \de^t_ke^s_l + \de^t_le^s_k)\tag1\\
[e^s\odot e^t, e^p_w] &= \de^t_w e^p\odot e^s + \de^s_w e^p\odot
e^t.\tag2
\endalign$$
We shall express the deformation tensor $\Ga$ by its values on
the generators, so we write $\Ga.(e_i\odot e_j) =:
\sum_{s,t}\Ga_{(st)(ij)}e^s\odot e^t$. Similarly to the above cases we
compute the deformation of the curvature.
$$
\aligned
[\Ga.(e_k&\odot e_l), e_i\odot e_j]  - [\Ga.(e_i\odot e_j), e_k\odot
e_l] = \\ &=\tsum_{s,t} \bigl( \Ga_{(st)(kl)}.[e^s\odot e^t,
e_i\odot e_j] - \Ga_{(st)(ij)}[e^s\odot e^t, e_k\odot
e_l]\bigr)\\
&=-\frac14\tsum_{s,t}\bigl(\Ga_{(st)(kl)}(\de^s_ie^t_j + \de^s_j e^t_i
+\\&\newquad
\de^t_i e^s_j + \de^t_j e^s_i) - \Ga_{(st)(ij)}(\de^s_ke^t_l +
\de^s_l e^t_k + \de^t_k e^s_l + \de^t_l e^s_k)\bigr)\\
&= \frac12\tsum_{p,w}\bigl(\de^w_l\Ga_{(kp)(ij)} +
\de^w_k\Ga_{(lp)(ij)} - \de^w_j\Ga_{(ip)(kl)} - \de^w_i\Ga_{(jp)(kl)}
\bigr)e^p_w
\endaligned\tag3$$
In order to get the deformation of the trace we compute	
$\partial^*(\ka_0)(e_k\odot e_l)$:
$$\aligned
[e^i\odot e^j,[\Ga.(e_k&\odot e_l), e_i\odot e_j]  -
[\Ga.(e_i\odot e_j), e_k\odot e_l]] = \\
&= \frac12\tsum_{p,w}(\text{the above coefficient at
$e^p_w$})(\de^j_we^p\odot e^i + \de^i_w e^p\odot e^j)\\
&=\tsum_{p,q}\bigl(\Ga_{(kp)(ql)} + \Ga_{(lp)(qk)} -
(m+1)\Ga_{(pq)(kl)}\bigr)e^p\odot e^q
\endaligned\tag4$$
and so the value of the deformation of the trace on the generators is
$$
\de\Tr(\ka_0)_{(pq)(kl)}=\de\Tr(\ka_0)(e_k\odot e_l, e_p\odot e_q) =
\Ga_{(kp)(ql)} + \Ga_{(lp)(qk)} - (m+1)\Ga_{(pq)(kl)}.
\tag5$$
Now, similarly to the Grassmannian case, we have to consider a
suitable combination. Surprisingly enough, we do not need the traces
inside of $\goth g_0$ in order to express the tensor $\Ga$.
If we substitute (5) into
$$m\de\Tr(\ka_0)_{(pq)(kl)} +
\de\Tr(\ka_0)_{(pk)(ql)} + \de\Tr(\ka_0)_{(pl)(qk)},$$
we are left
with $(2-m(m+1))\Ga_{(pq)(kl)}$ on the right hand side. Thus if we
start with a linear harmonic connection $\ga$ on $M$ and $\ka$ is the
curvature of $\tilde \ga$, then we can
achieve vanishing of $\partial^*\bar\ka_0$ on the section which
corresponds to $\ga$ by the choice
$$
\Ga_{(pq)(kl)} = \tfrac{1}{m(m+1)-2}(m \Tr(R)_{(pq)(kl)} +
\Tr(R)_{(pk)(ql)} + \Tr(R)_{(pl)(qk)})
\tag6$$
where $\Tr(R)$ is the Ricci curvature of $\ga$. In view of Lemma
\nmb!{3.3}, this deformation tensor leads to a normal Cartan
connection. Moreover, this deformation is uniquely determined by our
computation.
Thus, we have proved

\proclaim{\nmb.{3.11}. Theorem} Let $P\to M$ be an
almost Lagrangian structure
with a harmonic torsion, on a smooth manifold $M$ with dimension
greater then 2.
Then there is a uniquely defined normal Cartan connection $\om$
on $P$ and for each linear harmonic connection $\ga$ on the
underlying first order structure $P_0$ with curvature $R$,
$\om = \tilde\ga-\Ga\o\th_{-1}$, where the corresponding deformation tensor
$\Ga$ is given by the formula \nmb!{3.10}.(6).
\endproclaim

\subheading{\nmb.{3.12}. The almost spinorial case}
The situation is very similar to the almost Lagrangian case. We
have to proceed quite analogously with the symmetric matrices
replaced by the antisymmetric ones.

We have $\goth g_{-1}= \La^2\Bbb R^m$, $\goth g_1=\La^2\Bbb R^{m*}$,
$\goth g_0 = \goth g\goth l(m,\Bbb R)$.
The commutators of the base elements are
$$\align
[e^s\wedge e^t, e_k\wedge e_l] &= \tfrac14(-\de^s_le^t_k +
\de^t_le^s_k + \de^s_ke^t_l - \de^t_ke^s_l)\tag1\\
[e^s\wedge e^t, e^p_w] &= \de^t_w e^s\wedge e^p - \de^s_w e^t\wedge
e^p.\tag2
\endalign$$
We write, $\Ga.(e_i\wedge e_j) =:
\sum_{s,t}\Ga_{[st][ij]}e^s\wedge e^t$. Similarly as before we
compute the deformation of the curvature.
$$
\aligned
[\Ga.(e_k&\wedge e_l), e_i\wedge e_j]  - [\Ga.(e_i\wedge e_j), e_k\wedge
e_l] = \\ &=\tsum_{s,t} \bigl( \Ga_{[st][kl]}.[e^s\wedge e^t,
e_i\wedge e_j] - \Ga_{[st][ij]}[e^s\wedge e^t, e_k\wedge
e_l]\bigr)\\
&=\tfrac14\tsum_{s,t}\bigl(\Ga_{[st][kl]}(-\de^s_je^t_i + \de^t_j e^s_i
+ \de^s_i e^t_j - \de^t_i e^s_j) -\\&\newquad
\Ga_{[st][ij]}(-\de^s_le^t_k +
\de^t_l e^s_k + \de^s_k e^t_l - \de^t_k e^s_l)\bigr)\\
&= \tfrac12\tsum_{p,w}\bigl(\de^w_i\Ga_{[pj][kl]} +
\de^w_j\Ga_{[ip][kl]} + \de^w_k\Ga_{[lp][ij]} + \de^w_l\Ga_{[pk][ij]}
\bigr)e^p_w
\endaligned\tag3$$
Further we compute	
$$\aligned
\partial^*(\ka_0)(e_k\wedge e_l)&=[e^i\wedge e^j,[\Ga.(e_k\wedge e_l),
e_i\wedge e_j]  - [\Ga.(e_i\wedge e_j), e_k\wedge e_l]] = \\
&= \tfrac12\tsum_{p,w}(\text{the above coefficient at
$e^p_w$})(\de^j_we^i\wedge e_p - \de^i_w e^j\wedge e^p)\\
&=\tsum_{p,q}\bigl(\Ga_{[pk][ql]} + \Ga_{[lp][qk]} +
(1-m)\Ga_{[pq][kl]}\bigr)e^q\wedge e^p\\
\de\Tr(\ka_0)(e_k\wedge e_l&, e_p\wedge e_q) = \Ga_{[pk][ql]} +
\Ga_{[lp][qk]} - (m-1)\Ga_{[pq][kl]}
\endaligned\tag4$$
Now, we have to find a
suitable combination.
Let us substitute (4) into
$$m\de\Tr(\ka_0)_{[pq][kl]} +
\de\Tr(\ka_0)_{[pk][ql]} - \de\Tr(\ka_0)_{[pl][qk]}.$$
Then only
$(2-m(m-1))\Ga_{[pq][kl]}$ remains on the right hand side. Thus if we
start with a linear harmonic connection $\ga$ on $M$ and $\ka$ is the
curvature of $\tilde \ga$, then we can
achieve global vanishing of $\partial^*\bar\ka_0$ by the choice	
(cf\. Lemma \nmb!{3.3})
$$
\Ga_{[pq][kl]} = \tfrac{1}{m(m-1)-2}(m\Tr(R)_{[pq][kl]} +
\Tr(R)_{[pk][ql]} - \Tr(R)_{[pl][qk]})
\tag5$$
where $\Tr(R)$ is the Ricci curvature of $\ga$. Since
this deformation is uniquely determined by our
computation, we have proved:

\proclaim{\nmb.{3.13}. Theorem}
Let $P\to M$ be an
almost spinorial structure
with a harmonic torsion, on a smooth manifold $M$ with dimension
greater then 2.
Then there is a uniquely defined normal Cartan connection $\om$
on $P$ and for each linear harmonic connection $\ga$ on the
underlying first order structure $P_0$ with curvature $R$,
$\om = \tilde\ga-\Ga\o\th_{-1}$, where the corresponding deformation tensor
$\Ga$ is given by {\rm \nmb!{3.12}.(5)}.
\endproclaim

\enddocument